# Local existence of solutions and comparison principle for initial boundary value problem with nonlocal boundary condition for a nonlinear parabolic equation with memory

ALEXANDER GLADKOV

**Abstract:** We consider a nonlinear parabolic equation with memory $u_t = \Delta u + a u^p \int_0^t u^q(x,\tau)\,d\tau - b u^m$ for $(x,t) \in \Omega \times (0,+\infty)$ under nonlinear nonlocal boundary condition $\left.\dfrac{\partial u(x,t)}{\partial \nu}\right|_{\partial\Omega \times (0,+\infty)} = \int_\Omega k(x,y,t) u^l(y,t)\,dy$ and initial data $u(x,0) = u_0(x)$, $x \in \Omega$, where $a, b, q, m, l$ are positive constants, $p \geq 0$, $\Omega$ is a bounded domain in $\mathbb{R}^n$ with smooth boundary $\partial\Omega$, $\nu$ is unit outward normal on $\partial\Omega$. Nonnegative continuous function $k(x,y,t)$ is defined for $x \in \partial\Omega$, $y \in \bar\Omega$, $t \geq 0$, nonnegative function $u_0(x) \in C^1(\bar\Omega)$ and satisfies the condition $\dfrac{\partial u_0(x)}{\partial \nu} = \int_\Omega k(x,y,0) u_0^l(y)\,dy$ for $x \in \partial\Omega$. In this paper we study classical solutions. We establish the existence of a local maximal solution of the original problem. We introduce definitions of a supersolution and a subsolution. It is shown that under some conditions a supersolution is not less than a subsolution. We find conditions for the positiveness of solutions. As a consequence of the positiveness of solutions and the comparison principle of solutions, we prove the uniqueness theorem.

**Key words:** nonlinear parabolic equation; nonlocal boundary condition; existence of a solution; comparison principle

## Introduction

In this paper we consider the initial boundary value problem for the following nonlinear parabolic equation

$$u_t = \Delta u + au^p \int_0^t u^q(x,\tau)d\tau - bu^m, \quad x \in \Omega, \ t > 0, \qquad (1)$$

with nonlinear nonlocal boundary condition

$$\frac{\partial u(x,t)}{\partial \nu} = \int_\Omega k(x,y,t)u^l(y,t)dy, \quad x \in \partial\Omega, \ t > 0, \qquad (2)$$

and initial datum

$$u(x,0) = u_0(x), \quad x \in \Omega, \qquad (3)$$

where $p \geq 0, q > 0, m > 0, l > 0, a > 0, b > 0$, $\Omega$ is a bounded domain in $\mathbb{R}^n$ ($n \geq 1$) with smooth boundary $\partial\Omega$, $\nu$ is unit outward normal on $\partial\Omega$.

Throughout this paper we suppose that the functions $k(x,y,t)$ and $u_0(x)$ satisfy the following conditions:

$$k(x,y,t) \in C(\partial\Omega \times \bar\Omega \times [0,+\infty)), \ k(x,y,t) \geq 0;$$

$$u_0(x) \in C^1(\bar\Omega), \ u_0(x) \geq 0 \ in \ \Omega, \ \frac{\partial u_0(x)}{\partial \nu} = \int_\Omega k(x,y,0)u_0^l(y)dy \ on \ \partial\Omega.$$

Initial boundary value problems with nonlocal terms in parabolic equations or in boundary conditions have been considered in many papers (see, for example, [1–16] and the references therein). In particular, the initial boundary value problem (1)–(3) with $a = 0$ was considered for $b = b(x,t) \geq 0$ and $b = b(x,t) \leq 0$ in [17, 18] and [19, 20], respectively. The problem (1)–(3) with $p = 0$ and nonlocal boundary condition

$$u(x,t) = \int_\Omega k(x,y,t)u^l(y,t)dy, \quad x \in \partial\Omega, \ t > 0, \qquad (4)$$

was investigated in [21].

The aim of this paper is to study problem (1)–(3) for any positive $p, q, m$ and $l$. We prove existence of a local solution of (1)–(3). Comparison principle and the uniqueness of a solution are established. We show the nonuniqueness of solution of problem (1)–(3) with $u_0(x) \equiv 0$ also.

## Local existence

In this section a local existence theorem for (1)–(3) will be proved. We begin with definitions of a supersolution, a subsolution and a maximal solution of (1)–(3). Let $Q_T = \Omega \times (0,T), \ S_T = \partial\Omega \times (0,T), \ \Gamma_T = S_T \cup \bar\Omega \times \{0\}, \ T > 0$.

**Definition 1.** We say that a nonnegative function $u(x,t) \in C^{2,1}(Q_T) \cap C^{1,0}(Q_T \cup \Gamma_T)$ is a supersolution of (1)–(3) in $Q_T$, if

$$u_t \geq \Delta u + au^p \int_0^t u^q(x,\tau)d\tau - bu^m, \quad (x,t) \in Q_T, \tag{5}$$

$$\frac{\partial u(x,t)}{\partial \nu} \geq \int_\Omega k(x,y,t)u^l(y,t)dy, \quad (x,t) \in S_T, \tag{6}$$

$$u(x,0) \geq u_0(x), \quad x \in \Omega, \tag{7}$$

and $u(x,t) \in C^{2,1}(Q_T) \cap C^{1,0}(Q_T \cup \Gamma_T)$ is a subsolution of (1)–(3) in $Q_T$ if $u \geq 0$ and it satisfies (5)–(7) in the reverse order. We say that $u(x,t)$ is a solution of problem (1)–(3) in $Q_T$ if $u(x,t)$ is both a subsolution and a supersolution of (1)–(3) in $Q_T$.

**Definition 2.** We say that $u(x,t)$ is a maximal solution of (1)–(3) in $Q_T$ if for any other solution $w(x,t)$ of (1)–(3) in $Q_T$ the inequality $w(x,t) \leq u(x,t)$ is satisfied for $(x,t) \in Q_T \cup \Gamma_T$.

Let $\{\varepsilon_m\}$ be decreasing to 0 a sequence such that $0 < \varepsilon_m < 1$ and $\varepsilon_m \to 0$ as $m \to \infty$. For $\varepsilon = \varepsilon_m, m = 1,2,\ldots,$ let $u_{0\varepsilon}(x)$ be the functions with the following properties:

$$u_{0\varepsilon}(x) \in C^1(\overline{\Omega}), \; u_{0\varepsilon}(x) \geq \varepsilon, \; u_{0\varepsilon_i}(x) \geq u_{0\varepsilon_j}(x) \text{ for } \varepsilon_i \geq \varepsilon_j,$$
$$u_{0\varepsilon}(x) \to u_0(x) \text{ as } \varepsilon \to 0 \text{ uniformly in } \overline{\Omega}, \tag{8}$$
$$\frac{\partial u_{0\varepsilon}(x)}{\partial \nu} = \int_\Omega k(x,y,0)u_{0\varepsilon}^l(y)dy, \; x \in \partial\Omega.$$

Let us consider the following auxiliary problem:

$$\begin{cases} u_t = \Delta u + au^p \int_0^t u^q(x,\tau)d\tau - bu^m + b\varepsilon^m, & (x,t) \in Q_T, \\ \dfrac{\partial u(x,t)}{\partial \nu} = \int_\Omega k(x,y,t)u^l(y,t)dy, & (x,t) \in S_T, \\ u(x,0) = u_{0\varepsilon}(x), & x \in \Omega, \end{cases} \tag{9}$$

where $\varepsilon = \varepsilon_m$. The notion of a solution $u_\varepsilon$ for problem (9) can be defined in a similar way as in the Definition 1.

**Theorem 1.** Problem (9) has a unique solution in $Q_T$ for small values of $T > 0$.

*Proof.* Denote $K = \sup_{\partial\Omega \times Q_1} k(x,y,t)$ and introduce an auxiliary function $\psi(x)$ with the following properties:

$$\psi(x) \in C^2(\overline{\Omega}), \; \inf_\Omega \psi(x) \geq \max\left(\sup_\Omega u_{0\varepsilon}(x), 1\right), \; \inf_{\partial\Omega} \frac{\partial \psi(x)}{\partial \nu} \geq K\max(1,\exp(l-1))\int_\Omega \psi^l(y)dy.$$

We put

$$w(x,t) = \exp(\alpha t)\psi(x),$$

where $\alpha$ will be defined below.

To prove the existence of a solution for (9) we introduce the set

$$B = \{h(x,t) \in C(\bar{Q}_T) : \varepsilon \leq h(x,t) \leq w(x,t), h(x,0) = u_{0\varepsilon}(x)\}$$

and consider the following problem

$$\begin{cases} u_t = \Delta u + a\upsilon^p \int_0^t \upsilon^q(x,\tau)d\tau - bu^m + b\varepsilon^m, & (x,t) \in Q_T, \\ \dfrac{\partial u(x,t)}{\partial \nu} = \int_\Omega k(x,y,t)\upsilon^l(y,t)dy, & (x,t) \in S_T, \\ u(x,0) = u_{0\varepsilon}(x), & x \in \Omega, \end{cases} \quad (10)$$

where $\upsilon \in B$. It is obvious, $B$ is a nonempty convex subset of $C(\bar{Q}_T)$. By classical theory [22] problem (10) has a solution $u \in C^{2,1}(Q_T) \cap C^{1,0}(\bar{Q}_T)$ for small values of $T$. Let us call $A(\upsilon) = u$, where $\upsilon \in B$ and $u$ is a solution of (10). In order to show that $A$ has a fixed point in $B$ we verify that $A$ is a continuous mapping from $B$ into itself such that $AB$ is relatively compact. Obviously, the function $u(x,t) = \varepsilon$ is a subsolution of problem (10). Let us show that $w(x,t)$ is a supersolution of problem (10) for suitable choice of $\alpha > 0$ and $T > 0$.

Indeed,

$$w_t - \Delta w - a\upsilon^p \int_0^t \upsilon^q(x,\tau)d\tau + bw^m - b\varepsilon^m \geq w_t - \Delta w - aw^p \int_0^t w^q(x,\tau)d\tau + bw^m - b\varepsilon^m \geq$$

$$\geq \exp(\alpha t)[\alpha\psi(x) - \Delta\psi(x)] - a\exp(p\alpha t)\frac{\exp(q\alpha t) - 1}{q\alpha}\psi^{p+q} + b\left(\exp(m\alpha t)\psi^m(x) - \varepsilon^m\right) \geq 0$$

for $(x,t) \in Q_T$ if

$$\alpha \geq \max\left\{\frac{1}{q}, a\exp(1)\sup_\Omega \psi^{p+q-1}(x) + \sup_\Omega \frac{\Delta\psi(x)}{\psi(x)}\right\}, \quad T \leq \frac{1}{(p+q)\alpha}.$$

On the boundary $S_T$ we have for $T \leq 1/\alpha$

$$\frac{\partial w(x,t)}{\partial \nu} - \int_\Omega k(x,y,t)\upsilon^l(y,t)dy \geq \exp(\alpha t)K\max(1,\exp(l-1))\int_\Omega \psi^l(y)dy -$$

$$- K\exp(\alpha lt)\int_\Omega \psi^l(y)dy \geq 0.$$

The following inequality

$$w(x,0) - u_{0\varepsilon}(x) \geq 0$$

holds for $x \in \Omega$. Then $w(x,t)$ is a supersolution of problem (10) and thanks to a comparison principle for (10) $A$ maps $B$ into itself.

Let $G(x,y;t-\tau)$ denote the Green's function for a heat equation with homogeneous Neumann boundary condition. The Green's function has the following properties (see, for example, [23]):

$$G(x, y; t-\tau) \geq 0, \quad x, y \in \Omega, \ 0 \leq \tau < t,$$
$$\int_\Omega G(x, y; t-\tau) dy = 1, \quad x \in \Omega, \ 0 \leq \tau < t. \tag{11}$$

It is well known, that $u(x,t)$ is a solution of (10) in $Q_T$ if and only if for $(x,t) \in Q_T$

$$u(x,t) = \int_\Omega G(x,y;t) u_{0\varepsilon}(y) dy +$$
$$+ \int_0^t \int_\Omega G(x,y;t-\tau)\left(av^p(y,\tau)\int_0^\tau v^q(y,\sigma)d\sigma + b\left(\varepsilon^m - u^m(y,\tau)\right)\right) dy\, d\tau + \tag{12}$$
$$+ \int_0^t \int_{\partial\Omega} G(x,\xi;t-\tau) \int_\Omega k(\xi,y,\tau) v^l(y,\tau) dy\, dS_\xi\, d\tau.$$

We claim that $A$ is continuous. In fact let $v_k$ be a sequence in $B$ converging to $v \in B$ in $C(\bar{Q}_T)$. Denote $u_k = Av_k$. Then by (11), (12) we see that

$$|u - u_k| = \left|\int_0^t \int_\Omega G(x,y;t-\tau)\{a\left(v^p(y,\tau) - v_k^p(y,\tau)\right)\int_0^\tau v^q(y,\sigma)d\sigma + \right.$$
$$+ av_k^p(y,\tau)\int_0^\tau \left(v^q(y,\sigma) - v_k^q(y,\sigma)\right)d\sigma - b\left(u^m(y,\tau) - u_k^m(y,\tau)\right)\} dy\, d\tau +$$
$$+ \int_0^t \int_{\partial\Omega} G(x,\xi;t-\tau) \int_\Omega k(\xi,y,\tau)\left(v^l(y,\tau) - v_k^l(y,\tau)\right) dy\, dS_\xi\, d\tau \Big|$$
$$\leq aT^2 \sup_{Q_T} |v^p - v_k^p| \sup w^q + aT^2 \sup_{Q_T} |v^q - v_k^q| \sup w^p +$$
$$+ \theta T \sup_{Q_T} |u - u_k| + KT |\Omega| \sup_{Q_T} |v^l - v_k^l|,$$

where $\theta = mb \max\left(\varepsilon^{m-1}, \sup_{Q_T} w^{m-1}(x,t)\right)$, $T \leq \min\left\{1, \dfrac{1}{2\theta}\right\}$. Now we can conclude that $u_k$ converges to $u$ in $C(\bar{Q}_T)$ as $k \to \infty$.

The equicontinuity of $AB$ follows from (12) and the properties of the Green's function (see, for example, [24]). The Ascoli-Arzela theorem guarantees the relative compactness of $AB$. Thus we are able to apply the Schauder-Tychonoff fixed point theorem and conclude that $A$ has a fixed point in $B$ if $T$ is small. Now if $u_\varepsilon$ is a fixed point of $A$, $u_\varepsilon \in C^{2,1}(Q_T) \cap C^{1,0}(\bar{Q}_T)$ and it is a solution of (9) in $Q_T$. Uniqueness of the solution follows from a comparison principle for (9) which can be proved in a similar way as in the next section. The theorem is proved.

Now, let $\varepsilon_2 > \varepsilon_1$. Then it is easy to see that $u_{\varepsilon_2}(x,t)$ is a supersolution of problem (9) with $\varepsilon = \varepsilon_1$. Applying to problem (9) a comparison principle we have $u_{\varepsilon_1}(x,t) \leq u_{\varepsilon_2}(x,t)$. Using the last inequality and the continuation principle of solutions we deduce that the existence time of $u_\varepsilon$ does not decrease as $\varepsilon \to 0$. Taking $\varepsilon \to 0$, we get $u_M(x,t) = \lim_{\varepsilon \to 0} u_\varepsilon(x,t) \geq 0$ and $u_M(x,t)$ exists in $Q_T$ for some $T > 0$. We know that $u_\varepsilon(x,t)$ is a solution of (9) in $Q_T$ if and only if for $(x,t) \in Q_T$

$$u_\varepsilon(x,t) = \int_\Omega G(x,y;t)u_{0\varepsilon}(y)\,dy +$$
$$+ \int_0^t \int_\Omega G(x,y;t-\tau)\left(au_\varepsilon^p(y,\tau)\int_0^\tau u_\varepsilon^q(y,\sigma)\,d\sigma + b\left(\varepsilon^m - u_\varepsilon^m(y,\tau)\right)\right)dy\,d\tau + \quad (13)$$
$$+ \int_0^t \int_{\partial\Omega} G(x,\xi;t-\tau)\int_\Omega k(\xi,y,\tau)u_\varepsilon^l(y,\tau)\,dy\,dS_\xi\,d\tau.$$

Passing to the limit as $\varepsilon \to 0$ in (13), we obtain by dominated convergence theorem

$$u_M(x,t) = \int_\Omega G(x,y;t)u_0(y)\,dy +$$
$$+ \int_0^t \int_\Omega G(x,y;t-\tau)\left(au_M^p(y,\tau)\int_0^\tau u_M^q(y,\sigma)\,d\sigma - bu_M^m(y,\tau)\right)dy\,d\tau +$$
$$+ \int_0^t \int_{\partial\Omega} G(x,\xi;t-\tau)\int_\Omega k(\xi,y,\tau)u_M^l(y,\tau)\,dy\,dS_\xi\,d\tau$$

for $(x,t) \in Q_T$. Therefore, $u_M(x,t)$ is a solution of (1)–(3). Let $u(x,t)$ be any other solution of (1)–(3). Then by comparison principle from the next section $u_\varepsilon(x,t) \geq u(x,t)$. Taking $\varepsilon \to 0$, we conclude $u_M(x,t) \geq u(x,t)$. Now we proved the following local existence theorem.

**Theorem 2.** Problem (1)–(3) has a maximal solution in $Q_T$ for small values of $T$.

## Comparison principle

**Theorem 3.** Let $\bar{u}(x,t)$ and $\underline{u}(x,t)$ be a supersolution and a subsolution of problem (1)–(3) in $Q_T$, respectively. Suppose that $\underline{u}(x,t) > 0$ or $\bar{u}(x,t) > 0$ in $Q_T \cup \Gamma_T$ if either $\min(q,l) < 1$ or $0 < p < 1$. Then $\bar{u}(x,t) \geq \underline{u}(x,t)$ in $Q_T \cup \Gamma_T$.

*Proof.* Suppose that $\min(p,q,l) \geq 1$. Let $T_0 \in (0,T)$ and $u_{0\varepsilon}(x)$ have the same properties as in (8) but only $u_{0\varepsilon}(x) \to \underline{u}(x,0)$ as $\varepsilon \to 0$ uniformly in $\bar{\Omega}$. We can construct a solution $u_M(x,t)$ of (1)–(3) with $u_0(x) = \underline{u}(x,0)$ in the following way $u_M(x,t) = \lim_{\varepsilon \to 0} u_\varepsilon(x,t)$, where $u_\varepsilon(x,t)$ is a solution of (9). To establish the theorem we will show that

$$\underline{u}(x,t) \leq u_M(x,t) \leq \bar{u}(x,t), \ (x,t) \in \bar{Q}_{T_0}. \quad (14)$$

We prove the second inequality in (14) only since the proof of the first one is similar. Let $\varphi(x,\tau) \in C^{2,1}(\bar{Q}_{T_0})$ be a nonnegative function such that

$$\frac{\partial \varphi(x,t)}{\partial \nu} = 0 \ \text{for} \ (x,t) \in S_{T_0}.$$

If we multiply the first equation in (9) by $\varphi(x,t)$ and then integrate over $Q_t$ for $t \in (0,T_0)$, we obtain

$$\int_0^t \int_\Omega u_{\varepsilon\tau}(x,\tau)\varphi(x,\tau)\,dx\,d\tau =$$

$$= \int_0^t \int_\Omega \left( \Delta u_\varepsilon(x,\tau) + a u_\varepsilon^p(x,\tau) \int_0^\tau u_\varepsilon^q(x,\sigma)\,d\sigma + b\left(\varepsilon^m - u_\varepsilon^m(x,\tau)\right) \right)\varphi(x,\tau)\,dx\,d\tau.$$

Integrating by parts and using Green's identity, we have

$$\int_\Omega u_\varepsilon(x,t)\varphi(x,t)\,dx \leq \int_\Omega u_\varepsilon(x,0)\varphi(x,0)\,dx +$$
$$+ \int_0^t \int_\Omega \left( u_\varepsilon(x,\tau)\varphi_\tau(x,\tau) + u_\varepsilon(x,\tau)\Delta\varphi(x,\tau) \right) dx\,d\tau +$$
$$+ \int_0^t \int_\Omega \left( a u_\varepsilon^p(x,\tau) \int_0^\tau u_\varepsilon^q(x,\sigma)\,d\sigma + b\left(\varepsilon^m - u_\varepsilon^m(x,\tau)\right) \right)\varphi(x,\tau)\,dx\,d\tau + \qquad (15)$$
$$+ \int_0^t \int_{\partial\Omega} \varphi(x,\tau) \int_\Omega k(x,y,\tau) u_\varepsilon^l(y,\tau)\,dy\,dS_x\,d\tau.$$

On the other hand, $\bar{u}$ satisfies (15) with reversed inequality and with $\varepsilon = 0$. Set $w(x,t) = u_\varepsilon(x,t) - \bar{u}(x,t)$. Then $w(x,t)$ satisfies

$$\int_\Omega w(x,t)\varphi(x,t)\,dx \leq \int_\Omega w(x,0)\varphi(x,0)\,dx + \varepsilon^m b \int_0^t \int_\Omega \varphi(x,\tau)\,dx\,d\tau +$$
$$+ \int_0^t \int_\Omega w(x,\tau)\left( \varphi_\tau(x,\tau) + \Delta\varphi(x,\tau) - mb\theta_1^{m-1}(x,\tau) \right)\varphi(x,\tau)\,dx\,d\tau +$$
$$+ \int_0^t \int_\Omega \left( a\bar{u}^p(x,\tau)\varphi(x,\tau)\int_0^\tau q\theta_2^{q-1}(x,\sigma)w(x,\sigma)\,d\sigma \right) dx\,d\tau + \qquad (16)$$
$$+ \int_0^t \int_\Omega \left( ap\theta_3^{p-1}(x,\tau)w(x,\tau)\varphi(x,\tau) \int_0^\tau u_\varepsilon^q(x,\sigma)\,d\sigma \right) dx\,d\tau +$$
$$+ \int_0^t \int_{\partial\Omega} \varphi(x,\tau) \int_\Omega k(x,y,\tau) l\theta_4^{l-1}(y,\tau)w(y,\tau)\,dy\,dS_x\,d\tau,$$

where $\theta_i(x,\tau)$ ($i=1,2,3,4$) are some continuous functions between $u_\varepsilon(x,t)$ and $\bar{u}(x,t)$. Note here that by hypotheses for $k(x,y,t)$, $u_\varepsilon(x,t)$ and $\bar{u}(x,t)$, we have

$$0 \leq \bar{u}(x,t) \leq M,\ \varepsilon \leq u_\varepsilon(x,t) \leq M,\ (x,t) \in \bar{Q}_{T_0},$$
$$0 \leq k(x,y,t) \leq M,\ (x,y,t) \in \partial\Omega \times \bar{\Omega} \times [0,T_0], \qquad (17)$$

where $M$ is some positive constant. Then it is easy to see from (17) that $\theta_1^{m-1}(x,\tau)$, $\theta_2^{q-1}(x,\tau)$, $\theta_3^{p-1}(x,\tau)$ and $\theta_4^{l-1}(x,\tau)$ are positive and bounded functions in $\bar{Q}_{T_0}$ and, moreover, $\theta_1^{m-1}(x,\tau) \leq \max\{\varepsilon^{m-1}, M^{m-1}\}$, $\theta_2^{q-1}(x,\tau) \leq M^{q-1}$, $\theta_3^{p-1}(x,\tau) \leq M^{p-1}$, $\theta_4^{l-1}(x,\tau) \leq M^{l-1}$. Define a sequence $\{a_n\}$ in the following way: $a_n(x,t) \in C^\infty(\bar{Q}_{T_0})$, $a_n(x,t) \geq 0$ and $a_n(x,t) \to mb\theta_1^{m-1}(x,t)$ as $n \to \infty$ in $L^1(\bar{Q}_{T_0})$. Now, we consider a backward problem given by

$$\begin{cases} \varphi_\tau + \Delta\varphi - a_n\varphi = 0, \ (x,\tau) \in Q_t, \\ \dfrac{\partial \varphi(x,\tau)}{\partial \nu} = 0, \ (x,\tau) \in S_t, \\ \varphi(x,t) = \psi(x), \ x \in \Omega, \end{cases} \tag{18}$$

where $\psi(x) \in C_0^\infty(\Omega)$ and $0 \le \psi(x) \le 1$. Denote a solution of (18) as $\varphi_n(x,\tau)$. Then by the standard theory for linear parabolic equations (see [24], for example), we find that $\varphi_n(x,\tau) \in C^{2,1}(\bar{Q}_t)$, $0 \le \varphi_n(x,\tau) \le 1$ in $\bar{Q}_t$. Putting $\varphi = \varphi_n$ in (16) and passing to the limit as $n \to \infty$, we infer

$$\int_\Omega w(x,t)\psi(x)\,dx \le \int_\Omega w_+(x,0)\,dx + \varepsilon^m b T_0 |\Omega| + \\ + \{a(p+q)M^{p+q-1}T_0 + l|\partial\Omega|M^l\} \int_0^t \int_\Omega w_+(x,\tau)\,dx\,d\tau, \tag{19}$$

where $w_+ = \max(w,0)$, $|\partial\Omega|$ and $|\Omega|$ are the Lebesgue measures of $\partial\Omega$ in $\mathbb{R}^{n-1}$ and $\Omega$ in $\mathbb{R}^n$, respectively. Since (19) holds for every $\psi(x)$, we can choose a sequence $\{\psi_n(x)\}$ converging in $L^1(\Omega)$ to

$$\psi(x) = \begin{cases} 1, \text{ if } w(x,t) > 0, \\ 0, \text{ if } w(x,t) \le 0. \end{cases}$$

Passing to the limit as $n \to \infty$ in (19), we obtain

$$\int_\Omega w(x,t)_+\,dx \le \int_\Omega w_+(x,0)\,dx + \varepsilon^m b T_0 |\Omega| + \\ + \{a(p+q)M^{p+q-1}T_0 + l|\partial\Omega|M^l\} \int_0^t \int_\Omega w_+(x,\tau)\,dx\,d\tau, \ t \in (0,T_0].$$

Applying now Gronwall's inequality, we have for $t \in (0,T_0]$

$$\int_\Omega w_+(x,t)\,dx \le \left(\int_\Omega w_+(x,0)\,dx + \varepsilon^m b T_0 |\Omega|\right) \exp\left[\{a(p+q)M^{p+q-1}T_0 + l|\partial\Omega|M^l\}t\right].$$

Passing to the limit as $\varepsilon \to 0$, the conclusion of the theorem follows for $\min(p,q,l) \ge 1$. For the case $p = 0, \min(q,l) \ge 1$ we prove the theorem in the same way. If $\min(q,l) < 1$ or $0 < p < 1$ we can consider $w(x,t) = \underline{u}(x,t) - \bar{u}(x,t)$ and prove the theorem in a similar way using the positiveness of a subsolution or a supersolution. The theorem is proved.

**Remark 1.** For similar problem (1), (3), (4) with $p = 0$ the authors of [21] suppose in the comparison principle that $\underline{u}(x,t) > 0$ or $\bar{u}(x,t) > 0$ in $Q_T \cup \Gamma_T$ if $\min(q,m,l) < 1$.

**Lemma 1.** Let $u(x,t)$ be a solution of problem (1)–(3) in $Q_T$. Let $u_0(x) \not\equiv 0$ in $\Omega$ and $m \ge 1$. Then $u(x,t) > 0$ in $Q_T \cup S_T$. If $u_0(x) > 0$ in $\bar{\Omega}$ and $p < m < 1$ then $u(x,t) > 0$ in $Q_T \cup \Gamma_T$.

*Proof.* Let $u_0(x) \not\equiv 0$ in $\Omega$ and $m \ge 1$. We denote

$$M = \sup_{Q_{T_0}} u(x,t),$$

where $M$ is some positive constant, $T_0 \in (0,T)$. Now we put $h(x,t) = u(x,t)\exp(\lambda t)$ with $\lambda \geq bM^{m-1}$. Then in $Q_{T_0}$ we have

$$h_t - \Delta h = \exp(\lambda t)\left(\lambda u + u_t - \Delta u\right) \geq u\exp(\lambda t)\left(\lambda - bu^{m-1}\right) \geq 0.$$

Since $h(x,0) = u_0(x) \geq 0$, $x \in \Omega$, and $u_0(x) \not\equiv 0$ in $\Omega$, by the strong maximum principle $h(x,t) > 0$ in $Q_{T_0}$. Hence, $u(x,t) > 0$ in $Q_{T_0}$. Let $h(x_0,t_0) = 0$ in some point $(x_0,t_0) \in S_T$. Then according to Theorem 3.6 of [25] it yields $\partial h(x_0,t_0)/\partial \nu < 0$ which contradicts (2).

Let $u_0(x) > 0$ in $\overline{\Omega}$ and $p < m < 1$. Then there exist $\tau \in (0,T)$ and $\varepsilon > 0$ such that

$$u(x,t) \geq \varepsilon \text{ in } \overline{Q_\tau},$$

and, moreover, $u(x,t) \equiv \varepsilon_1 = \min(\varepsilon,[a\tau\varepsilon^q/b]^{1/(m-p)})$ is the subsolution of problem (1)–(3) in $Q_{T_0} \setminus \overline{Q_\tau}$ with initial function $u(x,\tau)$ for $t = \tau$ instead of (3). Putting $\underline{u}(x,t) \equiv \varepsilon_1$ and $\overline{u}(x,t) \equiv u(x,t)$ and arguing as in the proof of Theorem 3, we get

$$u(x,t) \geq \varepsilon_1 \text{ in } \overline{Q_{T_0}} \text{ for any } T_0 \in (0,T).$$

Lemma 1 is proved.

As a simple consequence of Theorem 3 and Lemma 1, we get the following uniqueness result for problem (1)–(3).

**Теорема 4.** Let problem (1)–(3) have a positive in $Q_T \cup \Gamma_T$ solution or a solution in $Q_T$ either with nonnegative initial data in $\Omega$ for $\min(p,q,l) \geq 1$ or with positive initial data in $\overline{\Omega}$ under the conditions $m \geq 1$ or $p < m < 1$. Then a solution of problem (1)–(3) is unique in $Q_T$.

Now we will prove the nonuniqueness of solution of problem (1)–(3) with $u_0(x) \equiv 0$ for $l < \min(1,m)$ or $p + q < \min(1,m)$. We note that problem (1)–(3) with $u_0(x) \equiv 0$ has solution $u(x,t) \equiv 0$.

**Теорема 5.** Let $u_0(x) \equiv 0$ and either $l < \min(1,m)$ and

$$k(x,y_0,t_0) > 0 \text{ for any } x \in \partial\Omega \text{ and some } y_0 \in \partial\Omega \text{ and } t_0 \in [0,T) \tag{20}$$

or $p + q < \min(1,m)$. Then maximal solution of problem (1)–(3) $u_M(x,t) \not\equiv 0$ in $Q_T$.

*Proof.* As shown in Theorem 2 a maximal solution $u_M(x,t) = \lim_{\varepsilon \to 0} u_\varepsilon(x,t)$, where $u_\varepsilon(x,t)$ is some positive in $\overline{Q_T}$ supersolution of (1)–(3). To prove the theorem we construct a subsolution $\underline{u}(x,t) \not\equiv 0$ of (1)–(3) with $u_0(x) \equiv 0$. By Theorem 3 we have $u_\varepsilon(x,t) \geq \underline{u}(x,t)$ and therefore maximal solution $u_M(x,t) \not\equiv 0$.

At first let $l < \min(1, m)$ and (20) holds. To construct a subsolution we use the change of variables in a neighborhood of $\partial \Omega$ as in [26]. Let $\bar{x}$ be a point on $\partial \Omega$. We denote by $\hat{n}(\bar{x})$ the inner unit normal to $\partial \Omega$ at the point $\bar{x}$. Since $\partial \Omega$ is smooth it is well-known that there exists $\delta > 0$ such that the mapping $\psi : \partial \Omega \times [0, \delta] \to \square^n$ given by $\psi(\bar{x}, s) = \bar{x} + s\hat{n}(\bar{x})$ defines new coordinates $(\bar{x}, s)$ in a neighborhood of $\partial \Omega$ in $\overline{\Omega}$.

Under the assumptions of the theorem, there exists $\bar{t}$ such that $k(x, y, t) > 0$ for $t_0 \leq t \leq t_0 + \bar{t}$, $x \in \partial \Omega$ and $y \in V(y_0)$, where $V(y_0)$ is some neighborhood of $y_0$ in $\overline{\Omega}$.

Let $1/(1-l) < \alpha \leq 1/(1-m)$ for $m < 1$ and $\alpha > 1/(1-l)$ for $m \geq 1$, $2 < \beta < 2/(1-m)$ for $m < 1$ and $\beta > 2$ for $m \geq 1$ and assume that $A > 0$, $0 < \xi_0 \leq 1$ and $0 < T_0 < \min(T - t_0, \bar{t}, \delta^2)$. For points in $\partial \Omega \times [0, \delta] \times (t_0, t_0 + T_0]$ of coordinates $(\bar{x}, s, t)$ define

$$\underline{u}(\bar{x}, s, t) = A(t - t_0)^\alpha \left( \xi_0 - \frac{s}{\sqrt{t - t_0}} \right)_+^\beta$$

and extend $\underline{u}$ as zero to the whole of $\overline{Q_\tau}$ with $\tau = t_0 + T_0$. Arguing as in [17] we prove that $\underline{u}$ is the subsolution of (1)–(3) with $u_0(x) \equiv 0$ in $Q_\tau$.

Now we suppose that $p + q < \min(1, m)$. Then it is easy to check that $\underline{u}(x, t) = t^\gamma$ is the subsolution of (1)–(3) with $u_0(x) \equiv 0$ in $Q_\tau$ for small values of $\tau$ if

$$\gamma > \max\left( \frac{2}{1-(p+q)}, \frac{1}{m-(p+q)} \right).$$

Theorem 5 is proved.